\documentclass[12pt, reqno]{amsart}
\allowdisplaybreaks
\usepackage{amssymb,amsthm,amsmath}
\usepackage[numbers,sort&compress]{natbib}

\title[]{Global existence of weak solutions to the 3D incompressible axisymmetric Euler equations without swirl}
\author{Quansen Jiu,\,Jitao Liu,\,Dongjuan Niu}
\address[Quansen Jiu]{School of Mathematical Sciences, Capital Normal University, Beijing, 100048, P. R. China.}
\email{jiuqs@cnu.edu.cn}
\address[Jitao Liu]{College of Applied Sciences, Beijing University of Technology, Beijing 100124, P. R. China.}
\email{jtliu@bjut.edu.cn,\,\,\,jtliumath@qq.com}
\address[Dongniu Niu]{School of Mathematical Sciences, Capital Normal University, Beijing, 100048, P. R. China.}
\email{djniu@cnu.edu.cn}


\theoremstyle{plain}
\newtheorem{corollary}{Corollary}[section]
\newtheorem{theorem}{Theorem}[section]
\newtheorem{lemma}{Lemma}[section]
\newtheorem{proposition}{Proposition}[section]

\theoremstyle{definition}
\newtheorem{definition}{Definition}[section]

\newtheorem{remark}{Remark}[section]

\let\f=\frac
\let\p=\partial
\def\R{\Bbb R}

\def\no{\noindent}
\def\endproof{\hphantom{MM}\hfill\llap{$\square$}\goodbreak}
\newcommand{\beq}{\begin{equation}}
\newcommand{\eeq}{\end{equation}}
\newcommand{\ben}{\begin{eqnarray}}
\newcommand{\een}{\end{eqnarray}}
\newcommand{\beno}{\begin{eqnarray*}}
\newcommand{\eeno}{\end{eqnarray*}}

\begin{document}

%%% ----------------------------------------------------------------------

\begin{abstract}
In this paper, we mainly investigate the tridimensional incompressible axisymmetric Euler equations without swirl in the whole space. Specifically, we prove the global existence of weak solutions if the swirl component of initial vorticity $w_0^\theta$ satisfies that $\f{w_0^\theta}r\in L^1\cap L^p(\R^3)$ for some $p>1$. To achieve this goal, we establish the $L_{\rm loc}^{2+\alpha}(\R^3)$ estimate of velocity fields for some $\alpha>0$, which is {\it innovative} to the best of our knowledge. Our result extends previous work in the literature.
\end{abstract}
%%% ----------------------------------------------------------------------
\maketitle
%%% ----------------------------------------------------------------------
%%%%%%%%%%%%%%%%%%%%%%%%%%%%%%%%%%%%%%%%%%%%%%%%%%%%%%%%%%%%%%%%%%%%%%%%%%
%INTRODUCTION%%%%%%%%%%%%%%%%%%%%%%%%%%%%%%%%%%%%%%%%%%%%%%%%%%%%%%%%%%%%%
%%%%%%%%%%%%%%%%%%%%%%%%%%%%%%%%%%%%%%%%%%%%%%%%%%%%%%%%%%%%%%%%%%%%%%%%%%

\vskip .3in
\section{Introduction and main results}\hspace*{\parindent}

In this paper, we are concerned with the three-dimensional incompressible Euler equations
\begin{align} \label{equ1} \left\{
\begin{aligned}
&\p_tu+u\cdot \nabla u=-\nabla p,\\
&\nabla\cdot u=0,
\end{aligned}
\right. \end{align}
in the whole space $\R^3$ with initial data $u(0,x)=u_0(x)$, where $u=(u_1, u_2, u_3)$
and $p=p(x,t)$ represent the velocity fields and pressure respectively.

\vskip .1in
The mathematical study to the incompressible Euler equations takes a long history with a large amount of associated literature. For two-dimensional case, Wolibner \cite{Wolibner} obtained the global well-posedness of smooth solutions in 1933. Then, this work was extended by Yudovich \cite{Yudovich1},  who proved the existence and uniqueness for a certain class of weak solutions if the initial vorticity $w_0$ lies in $L^1\cap L^{\infty}(\R^2)$. Later,  under the assumption of $w_0\in L^1\cap L^p(\R^2)$ for some $p>1$, Diperna and Majda showed that the weak solutions exist globally in \cite{DiPerna2}.  Furthermore, if $w_0$ is a finite Radon measure with one sign, there are also many works about the global existence of  weak solutions,  which can be  referred to \cite{Delort,Majda2,EM,Liujg2} for details.  However, the global existence of smooth solutions for 3D incompressible Euler equations with smooth initial data is still an important open problem, with a large literature.

\vskip .1in
From mathematical point of view, in two-dimensional case, the corresponding vorticity $w=\p_2u_1-\p_1u_2$ is a scalar fields and satisfies the following transport equation
$$\p_tw+u\cdot\nabla w=0,$$
which infers that its $L^p$ norm is conserved for all time. Nevertheless, for the three-dimensional case, $w$ becomes a vector fields and the {\it vortex stretching} term $w\cdot\nabla u$ appears in the equations of vorticity
$$\p_tw+u\cdot\nabla w=w\cdot\nabla u,$$
where $w=\nabla\times u$.
The presence of {\it vortex stretching} term brings more difficulties to prove the global regularity, which is the main reason causing this problem open. Therefore, many mathematicians explore the flows with certain geometrical assumptions, which shorten the gap between 2D and 3D flows. One typical situation is the 3D axisymmetric flows.

\vskip .1in
Whereas, even with this particular structure, it still open to exclude the singularity which occurs (if there is) only on the axis of symmetry (see \cite{CKN}), even for the Cauchy problem. But if the swirl component of velocity fields, $u_{\theta}$, is trivial, Ukhovskii, Yudovich \cite{UY} and Saint Raymond \cite{SR} proved that the weak solutions of 3D incompressible axisymmetric Euler equations are regular for all time. It should be noted that under this assumption, the corresponding vorticity quantity $\f{w_\theta}r$ is transported by a divergence free vector fields, which makes the problem more close to the 2D case.

\vskip .1in
However, for vortex sheet initial data, the problem of existence remains open, which is quite different from the 2D case. Afterwards, many mathematicians are committed to looking for {\it a little stronger} assumptions on the initial vorticity (comparing with the vortex sheet initial data), which suffice to guarantee the global existence of weak solutions. Inspired by recent progress in this direction for helically symmetric flows \cite{JLN}, we are interested in the answers for axisymmetric flows. There is a large literature devoted to axisymmetric Euler flows without swirl. In 1997, D. Chae and N. Kim proved the global existence of weak solutions under the assumptions that $\f{w_0^{\theta}}r\in L^{p}(\R^3)$ for some $p>6/5$ in \cite{ChaeKim}. Later, D. Chae and O. Y. Imanuvilov \cite{Chae2} obtained the similar result by assuming $u_0\in L^2(\R^3)$ and $|\f{w_0^{\theta}}r|[1+({\rm log}^{+}|\f{w_0^{\theta}}r|)^{\alpha}]\in L^1(\R^3)$ with $\alpha>1/2$. Recently, Jiu, Wu and Yang \cite{Jiu3} also obtained the existence result under the assumptions that $u_0\in L^2(\R^3)$ and $\f{w_0^\theta}{r}\in L^1\cap L^p(\R^3)\,\,({\rm for\,\,some\,\,}p>1)$ by using the method of viscous approximations. Regarding other related works, one can refer to \cite{JL,JL2,Jtliu,LMNP, Shirota,Gang,Danchin,Jiu2,Jiu1,LN,Anne,JLN,Titi,DiPerna3}.

\vskip .1in
It should be noted that in previous work, the initial assumption is not really a near-vortex-sheet data, because the initial velocity fields $u_0$ is assumed in $L^2(\R^3)$. The main reason lies in that their proofs are based on a key estimate
\ben\label{key1}
\int_0^T\int_{\R^3}\f1{1+z^2}\big(\f{u_r}{r}\big)^2dxdt\leq C\big(\|u_0\|_{L^2(\R^3)}^2+\|\f{w_0^{\theta}}r\|_{L^1(\R^3)}\big)
\een
raised by Chae-Imanuvilov \cite{Chae2}. As a matter of fact, for this model, whether or not $\f{w^{\theta}}r\in L^1\cap L^p(\R^3)\,(p>1)$ can imply $u\in L^2(\R^3)$, even $L_{\rm loc}^2(\R^3)$, is an interesting and open problem itself. Therefore, a natural question comes, given the initial vorticity $\f{w^{\theta}}r\in L^1\cap L^p(\R^3)\,(p>1)$, whether corresponding weak solutions exist globally?

\vskip .1in
In the present paper, we will give a positive answer to this question. That is, given initial data $\f{w_0^{\theta}}r\in L^1\cap L^{p}(\R^3)$ for some $p>1$, the weak solutions exist globally. What's more, there is a new and important observation that $\f{w^{\theta}}r\in L^1\cap L^p(\R^3)$ can imply $u\in L_{\rm loc}^{\f{2p}{2-p}}(\R^3)$ for $1<p<2$. From the point of view of mathematics, our work has been more close to the open problem investigated in \cite{Jiu1}, that is the so-called vortex sheet initial data problem for the 3D incompressible axisymmetric Euler equations without swirl.

\vskip .1in
In the process of proof, there are two big challenges in solving this problem without initial assumption $u_0\in L^2(\R^3)$. Firstly, the basic energy estimates take no effect and hence we do not have any estimates of velocity fields itself. To overcome them, we make the first attempt to establish the $L_{\rm loc}^p(\R^3)\,\,(p>1)$ estimate for the velocity fields, which is {\it new} to our knowledge. More precisely, we find out the explicit form of stream function in terms of vorticity. On the basis of this formulation, we then establish the estimate $\|u\|_{L^p_{\rm loc}(\R^3)}$ for any $p>1$. Then, it is natural to build up the $W_{\rm loc}^{1,p}(\R^3)\,(p>1)$ estimates of velocity fields. Yet, this is still far from resolving original problem, because current estimates only deduce the strong convergence of approximating solutions in $L^{2}(0,T; Q)$ for any $Q\subset\subset {\R^3}\backslash\{x\in R^3|r = 0\}$, other than
$L^2(0, T; L_{\rm loc}^2(\R^3))$.  As usual in previous work \cite{Jiu3}, current argument is enough to conclude the global existence of weak solutions, if the following propositon introduced by Jiu and Xin \cite{Jiu1} works.

{\bf Proposition.} {\it Suppose $u_0\in L^2(\R^3)$, for the approximate solutions $\{u^{\epsilon}\}$ constructed in Theorem \ref{exist}, if there exists a subsequence $\{u^{\epsilon_j}\} \subset\{u^{\epsilon}\}$ such that, for any $Q\subset\subset {\R^3}\backslash\{x\in R^3|r = 0\}$ and $\epsilon_j\rightarrow0,$
\beno
u^{\epsilon_j}\rightarrow u\quad {\rm \it strongly}\quad {\rm in}\quad L^2(0, T; L^2(Q)),
\eeno
then there exists a further subsequence of $\{u^{\epsilon_j}\}$, still denoted by itself, such
that, as $\epsilon_j\rightarrow0,$}
\beno
u^{\epsilon_j}\rightarrow u\quad {\rm \it strongly}\quad {\rm in}\quad L^2(0, T; L_{\rm loc}^2(\R^3)).
\eeno

\vskip .1in
However, in our case, this method would not work any more due to lack of initial assumption $u_0\in L^2(\R^3)$. This brings the other difficulty in solving this  problem, that is quite different from prior work. It then forces us to find a new way to establish the convergence of approximating solutions in the region contains the axis of symmetry. To solve it, we try to look for some  estimates of velocity fields {\it stronger} than $L^{2}_{\rm loc}(\R^3)$. In the end, we successfully established the $L^{\f{2p}{2-p}}_{\rm loc}(\R^3)$ estimates of velocity fields for $1<p<2$, based on delicate analysis of the structure of model and fully utilizing current a priori estimates. This estimate is {\it optimal} under current method and also a matter of concern. It should be mentioned that to prove it, we established the $L^{p}_{\rm loc}(\R_+^2)$ estimates of $\tilde{u}=(u_r,u_z)$, which is another {\it new} ingredient. Thus, on the basis of above facts, we finally establish the strong convergence of approximating solutions in $L^2(0, T; L_{\rm loc}^2(\R^3))$, which is sufficient to prove the global existence of weak solutions.

\vskip .1in
Before stating our main theorem, we would like to introduce the definition of weak solutions to (\ref{equ1}) as follow.

\begin{definition}[{\bf Global weak solutions}]\label{weakE} For any $T>0$, the velocity fields $u(x,t)\in L^\infty(0,T;L_{\rm loc}^{2}(\R^3))$ with initial data $u_0\in L_{\rm loc}^{1}(\R^3)$ is called a weak solution to the Euler equations (\ref{equ1}) if it holds\\
(i)\,\,For any vector field $\varphi\in C_0^\infty((0,T];\R^3)$ with $\nabla\cdot\varphi=0,$
\beno
&&\int_0^T\int_{\R^3}u\cdot\varphi_t+u\cdot\nabla \varphi\cdot u=\int_{\R^3}u_0\cdot\varphi_0\,dx,\notag
\eeno
(ii)\,\,For any function $\phi\in C_0^\infty((0,T];\R^3),$
\beno
&&\int_{\R^3}u\cdot\nabla\phi\,dx=0.
\eeno
\end{definition}

\vskip .1in
Under this definition, our main result can be summarized by the following theorem.

\begin{theorem}\label{thm41}
Suppose that $w_0^{\theta}=w_0^{\theta}(r,z)$ is a scalar axisymmetric function such that $w_0=w(x,0)=w_0^{\theta}e_{\theta}$ and $\f{w_0^\theta}r\in L^1\cap L^p(\R^3)$ for some $p>1$. Then, for any $T>0$, there exists at least an axisymmetric weak
solution $u$ without swirl in the sense of Definition \ref{weakE}.
\end{theorem}

This paper is organized as follows. In section 2, we introduce some notations and technical lemmas. In section 3, we will concentrate on the {\it a priori} estimates of velocity fields. Section 4 is devoted to proving the global existence of weak solutions, i.e., the proof of Theorem \ref{thm41}.

\vskip .3in
\section{Preliminary}\hspace*{\parindent}

In this section, we fix notations and set down some basic definitions. Initially, we would like to introduce the definition of 3D axisymmetric flow.

\begin{definition}[{\bf Axisymmetric flow}]\label{DAXI} A vector field $u(x,t)$ is called axisymmetric if it can be described by the form of
\ben\label{Axif}
u(x,t)=u_r(r,z,t)e_r+u_\theta(r,z,t)e_\theta+u_z(r,z,t)e_z
\een
in the cylindrical coordinate, where $e_r=(\hbox{cos}\theta,\hbox{sin}\theta,0)$, $e_\theta=(-\hbox{sin}\theta,\hbox{cos}\theta,0)$, $e_z=(0,0,1)
$. We call the velocity components $u_r(r,z,t),\,\,u_\theta(r,z,t),\,\,u_z(r,z,t)$ as radial, swirl and z-component respectively.
\end{definition}

{\bf Comment on notations:} In the following context, we will use $u_r,\,u_\theta,\,u_z$ to deonte $u_r(r,z,t),$ $u_\theta(r,z,t)$, $u_z(r,z,t)$ for simplicity.\\

With above definition, we set up the equations satisfied by $u_r,\,u_\theta,\,u_z$. Initially, under cylindrical coordinate, it is trivial that the gradient operator can be expressed in the form of ${\nabla}=e_r\p_r+\f1re_\theta \p_\theta+e_z\p_z$. Then, by some basic calculations, one can rewrite (\ref{equ1}) as
\begin{align} \label{equ2} \left\{
\begin{aligned}
&\p_tu_r+\tilde{u}\cdot\tilde{\nabla} u_r+\p_rp=\f{(u_\theta)^2}r,\\
&\p_tu_\theta+\tilde{u}\cdot\tilde{\nabla} u_\theta=-\f{u_\theta u_r}r,\\
&\p_tu_z+\tilde{u}\cdot\tilde{\nabla} u_z+\p_zp=0,\\
&\p_r(ru_r)+\p_z(ru_z)=0,
\end{aligned}
\right. \end{align}
where $\tilde{u}=(u_r, u_z)$ and $\tilde{\nabla}=(\p_r, \p_z)$. In addition, by $(\ref{equ2})^2$ and some basic calculations, it is clear that the quantity $ru_\theta$ satisfies the following transport equation:
\ben\label{equ3}
\p_t(r{u_\theta})+\tilde{u}\cdot\tilde{\nabla} (r{u_\theta})=0.
\een
Thanks to (\ref{equ3}), the following conclusion holds.

\begin{proposition} Assume $u$ is a smooth axisymmetric solution of 3D incompressible Euler equations, then the swirl component of velocity $u_\theta$ will be vanishing if its initial data $u_0^\theta$ be given zero.
\end{proposition}
\no{\bf Proof.}\quad Thanks to the incompressible condition $(\ref{equ2})^4$, by multiplying  (\ref{equ3}) with $r{u_\theta}$  and integrating  on $(0,t)$, it follows that
\beno
\|ru_\theta(t)\|_{L^2(\R^3)}\leq\|ru_0^\theta\|_{L^2(\R^3)}=0.
\eeno
Then, considering that $u_\theta$ is smooth and $u_\theta|_{r=0}\equiv 0$, we can conclude that $u_\theta\equiv0$ for any $t>0$.
\endproof

\vskip .1in
Therefore, if $u_0^\theta=0$, then the corresponding vector fields become $\tilde{u}$ and its voritcity can be described as $w=w_\theta e_\theta$, where $w_\theta=\p_zu_r-\p_ru_z$. What's more, the scalar quantity $\f{w_\theta}r$ is transported by $\tilde{u}$, i.e.,
\ben\label{equ4}
\p_t(\f{w_\theta}r)+\tilde{u}\cdot\tilde{\nabla} (\f{w_\theta}r)=0.
\een
This means that $\f{w_\theta}r$ is conserved along any particle trajectory. As a result, given the initial datum smooth sufficiently, the 3D axisymmetric Euler equations without swirl always possess a unique global solution (\cite{DiPerna1,SR}). Besides, by employing the incompressible condition and some basic calculations, we have the following conclusion.

\vskip .1in
{\bf Conservation laws for $\f{w_\theta}r$.} {\it Suppose $u$ is a smooth solution of 3D incompressible axisymmetric Euler equations, with its initial swirl component $u_0^\theta$ vanishing, then the estimate
\ben\label{est1}
\|\f{w_\theta}r\|_{L^p(\R^3)}\leq\|\f{w_0^\theta}r\|_{L^p(\R^3)}
\een
holds for any $p\in [1,\infty],$ where $w_0^\theta=w^\theta(x,0)$.\\
}

Subsequently, we will introduce the stream function, whose existence is proved in Lemma 2 of \cite{Liujg1}.

\begin{proposition}\label{stream1}
Let $u$ be a smooth axisymmetric velocity fields without swirl and $\nabla\cdot u=0$, then there exists a unique scalar function $\psi=\psi(r,z)$ such that
$u=\nabla\times(\psi e_\theta)$ and $\psi=0$ on the axis of symmetry $r=0.$
\end{proposition}

Finally, we will collect below some useful estimates of velocity fields in terms of $\f{w_\theta}{r}$, see \cite{Lei,JL,MZ} for instance.
\begin{lemma}\label{axi1}
Let $\psi$ be as in Proposition \ref{stream1}, it holds
\beno
\|\p_r^2\big(\f{\psi}r\big)\|_{L^p(\R^3)}+\|\f1r\p_r\big(\f{\psi}r\big)\|_{L^p(\R^3)}+\|\p_{rz}^2\big(\f{\psi}r\big)\|_{L^p(\R^3)}+\|\p_{z}^2\big(\f{\psi}r\big)\|_{L^p(\R^3)}\leq C\|\f{w_\theta}r\|_{L^p(\R^3)}
\eeno
for any $p>1$, where $C$ is an absolute constant. In particular,
\ben\label{p2}
\|\p_r\big(\f{u_r}r\big)\|_{L^p(\R^3)}+\|\p_z\big(\f{u_r}r\big)\|_{L^p(\R^3)}\leq C\|\f{w_\theta}r\|_{L^p(\R^3)}.
\een
\end{lemma}

\begin{lemma}\label{axi2}
Suppose that $u$ is a smooth solution of 3D incompressible axisymmetric Euler equations without swirl, then there holds
\ben\label{p3}
\|\f{u_r}r\|_{L^{\f{3p}{3-p}}}\leq C\|\f{w_\theta}r\|_{L^p}\quad\quad\quad\forall p\in(1,3),
\een
where $C$ is an absolute constant .
\end{lemma}

\vskip .3in
\section{A priori estimates of velocity fields}\hspace*{\parindent}

\subsection{$\mathbf {W_{\rm loc}^{1,p}(\R^3)\,(p>1)}$ estimates}\hspace*{\parindent}\\

In this section, we will focus on the $W_{\rm loc}^{1,p}(\R^3)$ estimates of velocity fields. Firstly, Proposition \ref{stream1} together with $\nabla\cdot u=0$ and $w=\nabla\times u=w_\theta e_\theta$, tells us that
$$-\Delta (\psi e_\theta)=w_\theta e_\theta.$$
Then by the elliptic theory, we have
\ben\label{stream21}
\psi(r_x,z_x) e_{\theta_x}=\int_{\R^3}G(X,Y)w_\theta(r_y,z_y) e_{\theta_y}dY,
\een
where $X=(r_x,\theta_x,z_x)$ and $G(X,Y)=|X-Y|^{-1}$ stands for the three-dimensional Green's function in the whole space. Regarding the Green's function $G(X,Y)$, it is well-known that the following two properties hold\\
{\bf (i)\,}:
\ben\label{stream22}
|D^k_X G(X,Y)|\leq C_k|X-Y|^{-1-k},
\een
{\bf (ii)\,}:
\ben\label{stream23}
G(\bar{X},Y)=G(X,\bar{Y}),\,\,\,\p_rG(\bar{X},Y)=\p_rG(X,\bar{Y}),\,\,\,\p_zG(\bar{X},Y)=\p_zG(X,\bar{Y})
\een
for all $(X,Y)\in\R^3$, $\bar{X}=(-x,-y, z)$ and $k=0, 1, 2$.

\vskip .1in
Until now, we have established the formulation (\ref{stream21}). However, in order to find out the explicit form of $\psi(r_x,z_x)$, we need to fix the value of $\theta_x$. Therefore, by making use of the rotational invariance and putting $\theta_x=0$ in (\ref{stream21}), we derive the explicit form of $\psi$ in terms of $w_{\theta}$
\ben\label{BSL}
\psi(r_x,z_x)=\int_{-\infty}^{\infty}\int_{0}^{\infty}\int_{-\pi}^{\pi}G(X,Y)w^\theta cos\theta_yr_yd{\theta_y}d{r_y}d{z_y},
\een
where $X=(r_x,0,z_x)$.

\vskip .1in
On this basis, we intend to utilize the stream function to establish $L_{\rm loc}^{p}(\R^3)$ estimates of velocity fields. And we would like to introduce the following Lemma, which is the cornerstone of this paper.

\begin{lemma}\label{stream2}
Assume $u$ and $\psi$ be as in Lemma \ref{stream1}, $w=\nabla\times u=w_\theta e_\theta$, then there holds that
\ben\label{stream200}
|\psi(r_x,z_x)|\leq C\int_{{\R}^3}{\rm min}\left(1,\frac{r_x}{|X-Y|}\right)\f{|w^\theta|}{|X-Y|}dY
\een
and
\ben\label{stream20}
|\p_r\psi(r_x,z_x)|+|\p_z\psi(r_x,z_x)|\leq C\int_{{\R}^3}{\rm min}\left(1,\frac{r_x}{|X-Y|}\right)\f{|w^\theta|}{|X-Y|^{2}}dY,
\een
where $C$ is an absolute constant and $X=(r_x,0,z_x)$.
\end{lemma}
\no{\bf Proof.}\quad
First of all, we do the estimate of $|\p_r\psi|.$ From (\ref{BSL}), we have
$$\p_r\psi=\int_{-\infty}^{\infty}\int_{0}^{\infty}\int_{-\pi}^{\pi}\p_rG(X,Y)w^\theta cos\theta_yr_yd{\theta_y}d{r_y}d{z_y},$$
which together with (\ref{stream23}) yields that
$$\p_r\psi=\int_{-\infty}^{\infty}\int_{0}^{\infty}\int_{-\f{\pi}2}^{\f{\pi}2}(\p_rG(X,Y)-\p_rG(\bar{X},Y))w^\theta cos\theta_yr_yd{\theta_y}d{r_y}d{z_y}.$$

Thus, to prove (\ref{stream20}), it suffices to verify that
\beno
&&H\triangleq\int_{-\f{\pi}2}^{\f{\pi}2}(\p_rG(X,Y)-\p_rG(\bar{X},Y))w^\theta cos\theta_yd{\theta_y}\\
&\leq& C\int_{-\f{\pi}2}^{\f{\pi}2}{\rm min}\left(1,\frac{r_x}{|X-Y|}\right)\f{|w^\theta|}{|X-Y|^{2}}d{\theta_y}.
\eeno
Without loss of generality, we assume $\theta^*$ be the unique real number $\theta_y\in[0,\f{\pi}2]$ such that $|X-Y|=r_x$ and split the integral $H$ into
$H=I+II+III$, with
$$I=\int_{-\f{\pi}2}^{-{\theta}^*} d{\theta_y},\quad II=\int_{-{\theta}^*}^{{\theta}^*}d{\theta_y},\quad III=\int_{{\theta}^*}^{\f{\pi}2}d{\theta_y},$$
where $|X-Y|>r_x$ for $I,\,III$ and $|X-Y|\leq r_x$ for $II$. Otherwise, $|X-Y|>r_x$ or $|X-Y|<r_x$ for all $\theta_y\in[-\f{\pi}2,\f{\pi}2]$. For these two cases, one can prove them along the same lines with estimating $I$ or $II$.

Because $|X-Y|\leq|\bar{X}-Y|$ for all $|\theta_y|\leq\f\pi2 $ and the interval $[-\theta^*, \theta^*]$ corresponds to those $\theta_y$ for which $|X-Y|\leq r_x$, one can conclude that $II$ satisfies the desired estimate easily.

Regarding the first and third terms, to start with, we fix some angle $\theta'\in [\theta^*, \f{\pi}2]$ and denote $X_{\beta}= (r{\rm cos}\,\beta,r{\rm sin}\,\beta,z)$ for
$\beta\in[-\pi,0]$. Besides, for the function $f(x,y,z)=f(r{\rm cos}\,\theta,r{\rm sin}\,\theta,z)$, it is clear that $\p_{\theta}f=r\p_h f\cdot e_\theta$, where ${\p_h=(\p_x,\p_y,0)}$. Therefore, by the fundamental theorem of calculus, it follows that
$$\p_rG(X,Y)-\p_rG(\bar{X},Y)=\pi r_x\int_{-{\pi}}^{0}\p_h\p_rG(X_{\beta},Y)\cdot e_{\beta} d{\beta}.$$
Then, by employing the fact $|X-Y|\leq|X_{\beta}-Y|$ for all ${\beta}\in[-\pi, 0]$ and (\ref{stream22}) , it holds that
$$|\p_rG(X,Y)-\p_rG(\bar{X},Y)|\leq Cr_x|X-Y|^{-3}.$$

Thus, we have obtained the estimate of $III$, that is
$$III\leq Cr_x\int_{\theta^*}^{\f{\pi}2}|X-Y|^{-3}|w^\theta|d{\theta_y}.$$
What's more, the estimate of $I$ can be treated by the same arguments with $III$. Thus, by adding up all the estimates, one can derive the estimate of $|\p_r\psi|$. As for $|\psi|$ and $|\p_z\psi|$, one can estimate it in the similar way and we will omit it here.
\endproof

\vskip .1in
Thanks to Lemma \ref{stream2}, we can then derive the upper bounds of $\f{\psi}r,\,\f{\p_r{\psi}}r,\,\f{\p_z{\psi}}r$ in terms of $\f{w_\theta}r$.

\begin{corollary}\label{stream3}
Under the assumptions of Lemma \ref{stream2}, it further holds that
\ben\label{stream31}
|\f{\psi(r_x,z_x)}{r_x}|\leq C\int_{{\R}^3}\f{|w^\theta|}{r_y|X-Y|}dY
\een
and
\ben\label{stream311}
|\f{\p_r\psi(r_x,z_x)}{r_x}|+|\f{\p_z\psi(r_x,z_x)}{r_x}|\leq C\int_{{\R}^3}\f{|w^\theta|}{r_y|X-Y|^2}dY,
\een
where $C$ is an absolute constant and $X=(r_x,0,z_x)$.
\end{corollary}
\no{\bf Proof.}\quad
Initially, if $Y\in\R^3$ are such that $|X-Y|\leq r_x$ for any $r_x$, then one has $r_y\leq r_x+|r_x-r_y|\leq r_x+|X-Y|\leq 2r_x$, which together with \eqref{stream200} and \eqref{stream20} implies that
\beno
|\f{\psi(r_x,z_x)}{r_x}|&\leq& C\int_{{\R}^3}{\frac{1}{r_x}}{\f{|w^\theta|}{|X-Y|}}dY\leq2C\int_{{\R}^3}{\frac{1}{r_y}}{\f{|w^\theta|}{|X-Y|}}dY
\eeno
and
\beno
|\f{\p_z\psi(r_x,z_x)}{r_x}|+|\f{\p_z\psi(r_x,z_x)}{r_x}|&\leq& C\int_{{\R}^3}{\frac{1}{r_x}}{\f{|w^\theta|}{|X-Y|^2}}dY\leq2C\int_{{\R}^3}{\frac{1}{r_y}}{\f{|w^\theta|}{|X-Y|^2}}dY.
\eeno
Otherwise, if $|X-Y|> r_x$, it is clear that $\f{r_y}{|X-Y|}\leq\f{r_x+|r_x-r_y|}{|X-Y|}\leq\f{r_x+|X-Y|}{|X-Y|}\leq 2$. Then, we can get that
\beno
|\f{\psi(r_x,z_x)}{r_x}|&\leq& C\int_{{\R}^3}\frac{1}{|X-Y|}{\f{|w^\theta|}{|X-Y|}}dY\leq2C\int_{{\R}^3}{\frac{1}{r_y}}{\f{|w^\theta|}{|X-Y|}}dY
\eeno
and
\beno
|\f{\p_z\psi(r_x,z_x)}{r_x}|+|\f{\p_z\psi(r_x,z_x)}{r_x}|&\leq& C\int_{{\R}^3}\frac{1}{|X-Y|}{\f{|w^\theta|}{|X-Y|^2}}dY\leq2C\int_{{\R}^3}{\frac{1}{r_y}}{\f{|w^\theta|}{|X-Y|^2}}dY.
\eeno
Thus, the proof is finished.
\endproof

\vskip .1in
\begin{remark}
The proof of Lemma \ref{stream2} and Corollary \ref{stream3} borrows some ideas from \cite{Shirota,Danchin}. In \cite{Danchin}, the authors used the explicit form of $|\f{\p_z\psi}r|$ in \eqref{stream311} to establish the $L^\infty(\R^3)$ estimate of $\f{u_r}r$. Here, we discover more applications of stream functions in doing some estimates of velocity fields, which will be shown in the following content.
\end{remark}

\vskip .1in
With the help of Lemma \ref{stream2} and Corollary \ref{stream3}, we can then derive the following {$L_{\rm loc}^{p}(\R^3)$ estimates} of velocity fields, which is the first key contribution of our work.

\begin{proposition}\label{estuu1}{\bf[$L_{\rm loc}^{p}(\R^3)$ estimates]}
Given $u$ as a smooth axisymmetric velocity fields without swirl satisfying $\nabla\cdot u=0$, then there holds
\beno
\|u\|_{L^p(B_{R}\times[-R,R])}
\leq C_{R}\|\f{w^{\theta}}r\|_{L^1\cap L^p(\R^3)}
\eeno
for any $p\in(1,\infty).$ Here $B_R=B_R(0)\subset\R^2$ be a 2D ball and the constant $C_R$ depends only on $R$.
\end{proposition}
\no{\bf Proof.}\quad According to Lemma \ref{stream1}, for axisymmetric smooth velocity fields $u$ with zero swirl component, there exists a unique stream function $\psi$ such that
$$u=u_re_r+u_ze_z=\nabla\times(\psi e_\theta).$$
This implies that $u_r=-\p_z\psi,\,\,u_z=\p_r\psi+\f{\psi}r$ and therefore $|u|\leq|\p_z\psi|+|\p_r\psi|+|\f{\psi}r|.$ Then by Lemma \ref{stream2} and Corollary \ref{stream3}, it follows that
\ben\label{estu12}
|u|&\leq&C\int_{{\R}^3}\f{|w^\theta|}{r_y|X-Y|}dY+C\int_{{\R}^3}\f{|w^\theta|}{|X-Y|^{2}}dY\notag\\
&\leq& C\int_{|X-Y|\leq1}\f{|w^\theta|}{r_y|X-Y|}dY+C\int_{|X-Y|>1}\f{|w^\theta|}{r_y|X-Y|}dY\notag\\
&&+C\int_{|X-Y|\leq1}\f{|w^\theta|}{|X-Y|^2}dY+C\int_{|X-Y|>1}\f{|w^\theta|}{|X-Y|^{2}}dY\notag\\
&\leq&C\int_{|X-Y|\leq1}\f{|w^\theta|}{r_y|X-Y|}dY+C\int_{|X-Y|>1}\f{|w^\theta|}{r_y|X-Y|}dY\notag\\
&&+Cr_x\int_{|X-Y|\leq1}\f{|w^\theta|}{r_y|X-Y|^{2}}dY+C\int_{|X-Y|\leq1}\f{|w^\theta||r_x-r_y|}{r_y|X-Y|^2}dY\notag\\
&&+Cr_x\int_{|X-Y|>1}\f{|w^\theta|}{r_y|X-Y|^{2}}dY+C\int_{|X-Y|>1}\f{|w^\theta||r_x-r_y|}{r_y|X-Y|^2}dY\notag\\
&\leq&C\int_{|X-Y|\leq1}\f{|w^\theta|}{r_y|X-Y|}dY+C\int_{|X-Y|>1}\f{|w^\theta|}{r_y|X-Y|}dY\\
&&+Cr_x\int_{|X-Y|\leq1}\f{|w^\theta|}{r_y|X-Y|^{2}}dY+C\int_{|X-Y|\leq1}\f{|w^\theta|}{r_y|X-Y|}dY\notag\\
&&+Cr_x\int_{|X-Y|>1}\f{|w^\theta|}{r_y|X-Y|^{2}}dY+C\int_{|X-Y|>1}\f{|w^\theta|}{r_y|X-Y|}dY\notag\\
&\leq&2C\int_{|X-Y|\leq1}\f{|w^\theta|}{r_y|X-Y|}dY+Cr_x\int_{|X-Y|\leq1}\f{|w^\theta|}{r_y|X-Y|^{2}}dY\notag\\
&&+2C\int_{|X-Y|>1}\f{|w^\theta|}{r_y|X-Y|}dY+Cr_x\int_{|X-Y|>1}\f{|w^\theta|}{r_y|X-Y|^{2}}dY\notag\\
&=&\sum\limits_{i=1}^4 I^i,\notag
\een
where we used the fact $|r_x-r_y|\leq |X-Y|$ in above inequalities. Therefore, by using of Young's inequality for convolutions, it holds that
\ben\label{estu13}
&&\|I^1\|_{L^p(B_{R}\times[-R,R])}+\|I^2\|_{L^p(B_{R}\times[-R,R])}\notag\\
&\leq&C\|\f{\chi_{\{|x|\leq1\}}}{|x|}\|_{L^1(\R^3)}\|\f{w^\theta}{r}\|_{L^p(\R^3)}+CR\|\f{\chi_{\{|x|\leq1\}}}{|x|^{2}}\|_{L^1(\R^3)}\|\f{w^\theta}{r}\|_{L^p(\R^3)}\\
&\leq&C(R+1)\|\f{w^\theta}{r}\|_{L^1\cap L^p(\R^3)}\notag
\een
for any $p\in(1,\infty)$ and cut-off function $\chi_{A}$ with compact support set $A$.

Regarding the left terms, by applying H\"{o}lder inequalities and Young's inequality for convolutions, it follows that
\ben\label{estu14}
&&\|I^3\|_{L^p(B_{R}\times[-R,R])}+\|I^4\|_{L^p(B_{R}\times[-R,R])}\notag\\
&\leq&CR^2\|I^3\|_{L^{{3p}}(B_{R}\times[-R,R])}+CR\|I^4\|_{L^{\f{3p}2}(B_{R}\times[-R,R])}\\
&\leq&CR^2\|\f{\chi_{\{|x|>1\}}}{|x|}\|_{L^{{3p}}(\R^3)}\|\f{w^\theta}{r}\|_{L^1(\R^3)}+CR^2\|\f{\chi_{\{|x|>1\}}}{|x|^{2}}\|_{L^{\f{3p}2}(\R^3)}\|\f{w^\theta}{r}\|_{L^1(\R^3)}\notag\\
&\leq&CR^2\|\f{w^\theta}{r}\|_{L^1(\R^3)}.\notag
\een
Finally, by summing up $(\ref{estu12})-(\ref{estu14})$, one can finish all the proof.
\endproof

\vskip .1in
Subsequently, we get to establish the {$L_{\rm loc}^{p}(\R^3)$ estimates} of $\nabla u$ in terms of $w_\theta.$ According to Proposition 2.20 in \cite{Majda1}, the gradient of velocity field can be expressed in terms of its vorticity by
\ben\label{singular1}
[\nabla u]h=[\mathcal{P}w]h+\f13w\times h.
\een
Here $\mathcal{P}$ is a singular integral operator of Calder\'{o}n-Zygmund type which is generated by a
homogeneous kernel of degree -3 (see \cite{Kato1}) and $h$ is a vector field. Moreover, the explicit form of $[\mathcal{P}w]h$ is
\ben\label{singular2}
[\mathcal{P}w]h=-P.V.\,\int_{\R^3}\big(\f1{4\pi}\f{w(y)\times h}{|x-y|^3}+\f3{4\pi}\f{\{[(x-y)\times w(y)]\otimes(x-y)\}h}{|x-y|^5}\big)dy.
\een

\vskip .1in
Therefore, with the help of (\ref{singular1}) and (\ref{singular2}), we are in the position to build up the following estimates.

\begin{proposition}\label{estuu2}{\bf[$\|\nabla u\|_{L^p_{\rm loc}(\R^3)}$ estimates]}
Assume that $u$ is a smooth axisymmterical velocity fields with divergence free and zero swirl component, then for any $p\in(1,\infty)$, there holds
\beno
\|\nabla u\|_{L^p(B_{R}\times[-R,R])}
\leq C_{R}\|\f{w^{\theta}}r\|_{L^1\cap L^p(\R^3)},
\eeno
where $B_R=B_R(0)\subset\R^2$ be a 2D ball and the constant $C_R$ depends only on $R$.
\end{proposition}
\no{\bf Proof.}\quad Thanks to (\ref{singular1}), it is clear that $\|\nabla u\|_{L^p}\simeq\sum\limits_{i}\|[\nabla u]e_i\|_{L^p}$ holds for any $p\in(1,\infty)$, where $e_{i} (i=r,\theta,z)$ is the orthogonal basis in (\ref{Axif}). Then, by setting $\chi(r,z)$ be a smooth cut-off function such that $\chi(r,z)=1$ in $B_{2R}\times[-2R,2R]$, and $\hbox{supp}\,\chi\subset B_{3R}\times[-3R,3R]$, we can split $[\nabla u]e_i$ into three parts as
\beno
[\nabla u]e_i&=&[\mathcal{P}(\chi w)]e_i+[\mathcal{P}\{(1-\chi)w\}]e_i+\f13w\times e_i\\
&=&I+II+III.
\eeno

Because $\mathcal{P}$ is a singular operator of Calder\'{o}n-Zygmund type, by the Calder\'{o}n-Zygmund inequality for $p\in(1,\infty)$, it is clear that
\ben\label{estu22}
&&\|I\|_{L^p(B_R\times[-R,R])}+\|III\|_{L^p(B_R\times[-R,R])}\notag\\
&\leq&C\|[\mathcal{P}(\chi w)]\|_{L^p(\R^3)}+C\|w\|_{L^p(B_R\times[-R,R])}\notag\\
&\leq&C\|w_\theta\|_{L^p(B_{2R}\times[-2R,2R])}\\
&\leq&CR\|\f{w_\theta}r\|_{L^p(\R^3)}.\notag
\een

As for the second term, by (\ref{singular2}), we have
\beno
II=-P.V.\,\int_{\R^3}\big(\f1{4\pi}\f{g(y)\times e_i}{|x-y|^3}+\f3{4\pi}\f{\{[(x-y)\times g(y)]\otimes(x-y)\}e_i}{|x-y|^5}\big)dy,
\eeno
where $g(y)=(1-\chi(y))w(y)$. In addition, as $\hbox{supp}\,(1-\chi(y))\subset \R^3\setminus B_{2R}\times[-2R,2R]$, it is clear that $|x-y|\geq |y|-|x|\geq R$ for $x\in B_{R}\times[-R,R]$ and $y\in\R^3\setminus B_{2R}\times[-2R,2R].$ Therefore, for $x\in B_{R}\times[-R,R]$, there holds
\beno
|II|&\leq&C\int_{|x-y|\geq R}\f{|w_\theta(y)|}{|x-y|^3}dy\\
&\leq&Cr_x\int_{|x-y|\geq R}\f{|w_\theta(y)|}{r_y|x-y|^3}dy+C\int_{|x-y|\geq R}\f{|w_\theta(y)||r_x-r_y|}{r_y|x-y|^3}dy\\
&\leq&Cr_x\int_{|x-y|\geq R}\f{|w_\theta(y)|}{r_y|x-y|^3}dy+C\int_{|x-y|\geq R}\f{|w_\theta(y)|}{r_y|x-y|^2}dy\\
&\leq&\f{C}{R^2}\|\f{w_\theta}r\|_{L^1(\R^3)},
\eeno
which further implies, after utilizing some basic calculations, that
\ben\label{estu23}
\|II\|_{L^p(B_R\times[-R,R])}\leq CR\|\f{w_\theta}r\|_{L^1(\R^3)}.
\een
Thus, we can finish the proof by adding up (\ref{estu22}) and (\ref{estu23}).
\endproof
\bigskip
\subsection{$\mathbf {L_{\rm loc}^{p}(\R^3)\,(p>2)}$ estimates}\hspace*{\parindent}\\

As stated in the introduction, to prove the global existence of weak solutions, we need the strong convergence of approximate solutions in $L^2(0, T; L_{\rm loc}^2(\R^3))$. Although we have built up the $W_{\rm loc}^{1,p}(\R^3)\,(p>1)$ estimates of velocity fields, it only imply the strong convergence of approximate solutions in $L^{2}(0,T; Q)$ for any $Q\subset\subset {\R^3}\backslash\{x\in R^3|r = 0\}$, other than $L^2(0, T; L_{\rm loc}^2(\R^3))$.

\vskip .1in
To solve this gap, we will focus on establishing the estimates of velocity fields {\it stronger} than $L^{2}_{\rm loc}(\R^3)$. The first step is to achieve the $L_{\rm loc}^{p}(\R_+^2)\,(p>1)$ estimates for $\tilde{u}$, which is a {\it new} ingredient in this manuscript.
\begin{lemma}\label{estuu3}{\bf[$\|\tilde{u}\|_{L^p_{\rm loc}(\R_+^2)}$ estimates]} Suppose $u=u_r(r,z,t)e_r+u_z(r,z,t)e_z$ is a smooth axisymmetric velocity fields without swirl satisfying $\nabla\cdot u=0$ and let $\tilde{u}=(u_r, u_z)$ , then the estimate
\beno
\|\tilde{u}\|_{L^p([0,R]\times[-R,R])}
\leq C_{R}\|\f{w^{\theta}}r\|_{L^1\cap L^p(\R^3)}
\eeno
holds for any $p\in(1,\infty)$ and the constant $C_R$ depending only on $R$.
\end{lemma}
\no{\bf Proof.}\quad Firstly, with the help of the estimate of $\|\f{u_r}r\|_{L^p(B_{R}\times[-R,R])}$ in Proposition \ref{estuu2} and noticing $p>1$, it is clear that
\ben\label{estu31}
&&\|u_r\|_{L^p([0,R]\times[-R,R])}\notag\\
&=&\big[\f1{2\pi}\int_{-R}^{R}\int_{0}^{R}\int_{-\pi}^{\pi}|\f{u_r}r|^pr^{p-1}rd{\theta}d{r}d{z}\big]^{\f1{p}}\notag\\
&\leq&CR^{1-\f1{p}}\|\f{u_r}r\|_{L^p(B_{R}\times[-R,R])}\\
&\leq&CR^{1-\f1{p}}\|\f{w^\theta}{r}\|_{L^1\cap L^p(\R^3)}\notag
\een

Regarding the estimate of $\|u_z\|_{L^p([0,R]\times[-R,R])}$, by Proposition \ref{stream1}, there holds that
$|u_z|\leq|\p_r\psi|+|\f{\psi}r|.$ Then we will estimate the two terms by different ways. For the first term, by similar skills as in \eqref{estu31} and Corollary \ref{stream3}, it follows that
$$\|{\p_r\psi}\|_{L^p([0,R]\times[-R,R])}\leq CR^{1-\f1p}\|\f{\p_r\psi}r\|_{L^p(B_{R}\times[-R,R])}$$
and
\ben\label{estu32}
&&|\f{\p_r\psi}r|\leq C\int_{{\R}^3}\f{|w^\theta|}{r_y|X-Y|^2}dY\notag\\
&\leq& C\int_{|X-Y|\leq1}\f{|w^\theta|}{r_y|X-Y|^2}dY+C\int_{|X-Y|>1}\f{|w^\theta|}{r_y|X-Y|^2}dY\\
&=&I_1+I_2.\notag
\een
Then, by making use of Young's inequality for convolutions, we finally deduce that
\ben\label{estu33}
&&\|{\p_r\psi}\|_{L^p([0,R]\times[-R,R])}\notag\\
&\leq&CR^{1-\f1{p}}\|I_1\|_{L^p(B_{R}\times[-R,R])}+CR\|I_2\|_{L^{\f{3p}2}(B_{R}\times[-R,R])}\\
&\leq&CR^{1-\f1{p}}\|\f{\chi_{\{|x|\leq1\}}}{|x|^2}\|_{L^1(\R^3)}\|\f{w^\theta}{r}\|_{L^p(\R^3)}+CR\|\f{\chi_{\{|x|>1\}}}{|x|^{2}}\|_{L^{\f{3p}2}(\R^3)}\|\f{w^\theta}{r}\|_{L^1(\R^3)}\notag\\
&\leq&C(R+1)\|\f{w^\theta}{r}\|_{L^1\cap L^p(\R^3)}\notag
\een
for any $p\in(1,\infty)$ and cut-off function $\chi_{A}$ with compact support set $A$. As for the other term, by using the notation $\tilde{X}=(r_x, z_x)$ and Corollary \ref{stream3}, we firstly obtain
\ben\label{estu34}
&&|\f{\psi}{r}|\leq C\int_{{\R}^3}\f{|w^\theta|}{r_y|X-Y|}dY\notag\\
&=&C\int_{-\infty}^{\infty}\int_{0}^{\infty}\int_{-\pi}^{\pi}\f{|w^\theta|}{\sqrt{r_x^2+r_y^2-2r_x r_y cos{\theta_y}+(z_x-z_y)^2}}d{r_y}d{\theta}d{z_y}\notag\\
&\leq&2\pi C\int_{-\infty}^{\infty}\int_{0}^{\infty}\f{|w^\theta|}{\sqrt{(r_x-r_y)^2+(z_x-z_y)^2}}d{r_y}d{z_y}\notag\\
&=&2\pi C\int_{-\infty}^{\infty}\int_{0}^{\infty}\f{|w^\theta|}{|\tilde{X}-\tilde{Y}|}d{r_y}d{z_y}\\
&\leq&C\int_{-\infty}^{\infty}\int_{0}^{2R}\f{|w^\theta|}{|\tilde{X}-\tilde{Y}|}d{r_y}d{z_y}+C\int_{-\infty}^{\infty}\int_{2R}^{\infty}\f{|w^\theta|}{|\tilde{X}-\tilde{Y}|}d{r_y}d{z_y}\notag\\
&\leq&C\int_{\R^2}\f{|w^\theta|\chi_{\{0<r_y<2R\}}}{|\tilde{X}-\tilde{Y}|}d\tilde{Y}+C\int_{-\infty}^{\infty}\int_{2R}^{\infty}\f{|w^\theta|}{|\tilde{X}-\tilde{Y}|}d{r_y}d{z_y}\notag\\
&\leq&C\int_{|\tilde{X}-\tilde{Y}|\leq 1}\f{|w^\theta|\chi_{\{0\leq r_y<2R\}}}{|\tilde{X}-\tilde{Y}|}d\tilde{Y}+C\int_{|\tilde{X}-\tilde{Y}|> 1}\f{|w^\theta|\chi_{\{0\leq r_y<2R\}}}{|\tilde{X}-\tilde{Y}|}d\tilde{Y}\notag\\
&&+C\int_{-\infty}^{\infty}\int_{2R}^{\infty}\f{|w^\theta|}{|\tilde{X}-\tilde{Y}|}d{r_y}d{z_y}\notag\\
&=&I_3+I_4+I_5,\notag
\een
where we used the fact that $w_{\theta}=0$ on the axis of symmetry $r=0$ in the fourth inequality. Then, for any $0\leq r_x<R$ and $r_y>2R$, it clear holds $|\tilde{X}-\tilde{Y}|>R$ and then $I_5\leq\f{C}R\|\f{w^\theta}{r}\|_{L^1(\R^3)}$. Thus, by applying Young's inequality for convolutions, we have
\ben\label{estu35}
&&\|\f{\psi}r\|_{L^p([0,R]\times[-R,R])}\notag\\
&\leq&C\|I_3\|_{L^p([0,R]\times[-R,R])}+CR^{\f1p}\|I_4\|_{L^{2p}([0,R]\times[-R,R])}+C\|I_5\|_{L^p([0,R]\times[-R,R])}\notag\\
&\leq&C\|I_3\|_{L^p(\R^2)}+CR^{\f1{p}}\|I_4\|_{L^{2p}(\R^2)}+CR\|\f{w^\theta}{r}\|_{L^1(\R^3)}\notag\\
&\leq&C\|\f{\chi_{\{|x|\leq1\}}}{|x|}\|_{L^1(\R^2)}\|{w^\theta}\chi_{\{0<r<2R\}}\|_{L^p(\R^2)}+CR^{\f1{p}}\|\f{\chi_{\{|x|>1\}}}{|x|}\|_{L^{2p}(\R^2)}\|{w^\theta}\chi_{\{0<r<2R\}}\|_{L^1(\R^2)}\notag\\
&&+CR\|\f{w^\theta}{r}\|_{L^1(\R^3)}\\
&\leq&CR^{1-\f1{p}}\|\f{w^\theta}{r}\|_{L^p(\R^3)}+CR^{\f1{p}}\|\f{w^\theta}{r}\|_{L^1(\R^3)}+CR\|\f{w^\theta}{r}\|_{L^1(\R^3)}\notag\\
&\leq&C(R+1)\|\f{w^\theta}{r}\|_{L^1\cap L^p(\R^3)},\notag
\een
which together with \eqref{estu33} further implies
\ben\label{estu36}
\|u_z\|_{L^p([0,R]\times[-R,R])}\leq C(R+1)\|\f{w^\theta}{r}\|_{L^1\cap L^p(\R^3)}.
\een
In the end, we can finish all the proof by adding up \eqref{estu31} and \eqref{estu36}.
\endproof

\vskip .1in
Thanks to Lemma \ref{estuu4} and by fully exploiting the structure of axisymmetric flows without swirl, we then build up the following estimates stronger than ${L^{2}_{\rm loc}(\R^3)}$.
\begin{proposition}\label{estuu4}{\bf[$\|{u}\|_{L^{\f{2p}{2-p}}_{\rm loc}(\R^3)}$ estimates]} Let $u$ be a smooth axisymmetric velocity fields without swirl as in Lemma \ref{estuu3}, then the estimate
\beno
\|{u}\|_{L^{\f{2p}{2-p}}(B_{R}\times[-R,R])}
\leq C_{R}\|\f{w^{\theta}}r\|_{L^1\cap L^p(\R^3)}
\eeno
holds for any $1<p<2$. Here $B_R=B_R(0)\subset\R^2$ is a 2D ball and the constant $C_R$ depending only on $R$.
\end{proposition}
\no{\bf Proof.}\quad {{\bf Step 1:} $\mathbf {u_r\in L^{\f{2p}{2-p}}_{\rm loc}(\R^3)}$}\quad Thanks to the Sobolev embedding inequality $W_{\rm loc}^{1,p}(\R^2_{+})$
$\hookrightarrow L_{\rm loc}^{\f{2p}{2-p}}(\R^2_{+})$ for any $1<p<2$, and the equality that
$$\|r^{\f{2-p}{2p}}u_r\|_{L^{\f{2p}{2-p}}([0,R]\times[-R,R])}={2\pi}^{-\f{2-p}{2p}}\|u_r\|_{L^{\f{2p}{2-p}}(B_R\times[-R,R])},$$
to prove $u_r\in{L^{\f{2p}{2-p}}_{\rm loc}(\R^3)}$, it suffices to verify $r^{\f{2-p}{2p}}u_r\in{W^{1,p}_{\rm loc}(\R_{+}^2)}$. First of all, we certify $r^{\f{2-p}{2p}}u_r\in L^p([0,R]\times[-R,R])$. Through some basic calculations and Proposition \ref{estuu2}, it clearly follows that
\ben\label{estu41}
&&\|r^{\f{2-p}{2p}}u_r\|_{L^p([0,R]\times[-R,R])}=\big[\f1{2\pi}\int_{-R}^{R}\int_{0}^{R}\int_{-\pi}^{\pi}|\f{u_r}r|^pr^{\f{p}2}rd{\theta}d{r}d{z}\big]^{\f1{p}}\notag\\
&\leq&CR^{\f12}\|\f{u_r}r\|_{L^p(B_{R}\times[-R,R])}\\
&\leq&CR^{\f12}\|\f{w^\theta}{r}\|_{L^1\cap L^p(\R^3)}.\notag
\een
In the second stage, we demonstrate $\p_r\big(r^{\f{2-p}{2p}}u_r\big)\in L^p([0,R]\times[-R,R])$. To achieve this goal, we decompose it into two terms by $\p_r\big(r^{\f{2-p}{2p}}u_r\big)=\p_r\big(\f{u_r}rr^{\f{2+p}{2p}}\big)=\p_r\big(\f{u_r}r\big)r^{\f{2+p}{2p}}+\f{2+p}{2p}\big(\f{u_r}r\big)r^{\f{2-p}{2p}}$ and estimate them separately. Again by some basic calculations and borrowing \eqref{p2} in Lemma \ref{axi1}, we have
\ben\label{estu42}
&&\|r^{\f{2+p}{2p}}\p_r\big(\f{u_r}r\big)\|_{L^p([0,R]\times[-R,R])}=\big[\f1{2\pi}\int_{-R}^{R}\int_{0}^{R}\int_{-\pi}^{\pi}|\p_r\big(\f{u_r}r\big)|^pr^{\f{p}2}rd{\theta}d{r}d{z}\big]^{\f1{p}}\notag\\
&\leq&CR^{\f12}\|\p_r\big(\f{u_r}r\big)\|_{L^p(B_{R}\times[-R,R])}\\
&\leq&CR^{\f12}\|\f{w^\theta}{r}\|_{L^1\cap L^p(\R^3)}.\notag
\een
The other term can be estimated by H\"{o}lder inequality and Lemma \ref{axi2}, that is
\ben\label{estu43}
&&\|{\f{2+p}{2p}}\big(\f{u_r}r\big)r^{\f{2-p}{2p}}\|_{L^p([0,R]\times[-R,R])}\leq\big[\f1{2\pi}\int_{-R}^{R}\int_{0}^{R}\int_{-\pi}^{\pi}|\f{u_r}r|^pr^{-\f{p}2}rd{\theta}d{r}d{z}\big]^{\f1{p}}\notag\\
&\leq&\big[\f1{2\pi}\int_{-R}^{R}\int_{0}^{R}\int_{-\pi}^{\pi}|\f{u_r}r|^{\f{3p}{3-p}}rd{\theta}d{r}d{z}\big]^{\f{3-p}{3}}\big[\f1{2\pi}\int_{-R}^{R}\int_{0}^{R}\int_{-\pi}^{\pi}r^{-\f{3}2}rd{\theta}d{r}d{z}\big]^{\f{1}{3}}\notag\\
&\leq&CR^{\f13}\|\f{u_r}r\|_{L^{\f{3p}{3-p}}(\R^3)}\big[\int_{0}^{R}r^{-\f{1}2}d{r}\big]^{\f{1}{3}}\\
&\leq&CR^{\f12}\|\f{w^\theta}{r}\|_{L^1\cap L^p(\R^3)}.\notag
\een
Regarding the term $\p_z\big(r^{\f{2-p}{2p}}u_r\big)$, due to $\p_z\big(r^{\f{2-p}{2p}}u_r\big)=\p_z\big(\f{u_r}r\big)r^{\f{2+p}{2p}}$, the way to estimate it would be along the same line with $\p_r\big(\f{u_r}r\big)r^{\f{2+p}{2p}}$ in \eqref{estu42} and we will omit it here to avoid repetition.\\
{{\bf Step 2:} $\mathbf {u_z\in L^{\f{2p}{2-p}}_{\rm loc}(\R^3)}$}\quad Through recalling Proposition \ref{stream1}, it is clear that
\ben\label{uz}
u_z=\p_r\psi+\f{\psi}{r}=r\p_r\big(\f{\psi}{r}\big)+\f{2\psi}{r}
\een
and we will deal with the two terms by different methods. For the term $\f{\psi}{r}$, we will estimate it by straightforward calculations.
According to Corollary \ref{stream3}, it yields
\ben\label{estu44}
&&|\f{\psi}r|\leq C\int_{{\R}^3}\f{|w^\theta|}{r_y|X-Y|}dY\notag\\
&\leq& C\int_{|X-Y|\leq1}\f{|w^\theta|}{r_y|X-Y|}dY+C\int_{|X-Y|>1}\f{|w^\theta|}{r_y|X-Y|}dY\\
&=&I_1+I_2,\notag
\een
which further implies, after making use H\"{o}lder inequality in bounded domain $B_{R}\times[-R,R]$ and Young's inequality for convolutions, that
\ben\label{estu45}
&&\|\f{\psi}r\|_{L^{\f{2p}{2-p}}(B_{R}\times[-R,R])}\leq C\|I_1\|_{L^{\f{2p}{2-p}}(\R^3)}+CR^{\f{6-3p}{4p}}\|I_2\|_{L^{\f{4p}{2-p}}(\R^3)}\notag\\
&\leq&C\|\f{\chi_{\{|x|\leq1\}}}{|x|}\|_{L^{2}(\R^3)}\|\f{w^\theta}{r}\|_{L^p(\R^3)}+C(1+R)\|\f{\chi_{\{|x|>1\}}}{|x|}\|_{L^{\f{4p}{2-p}}(\R^3)}\|\f{w^\theta}{r}\|_{L^1(\R^3)}\notag\\
&\leq&C(R+1)\|\f{w^\theta}{r}\|_{L^1\cap L^p(\R^3)}
\een
for $1<p<2$. In above inequalities, we have used $\f14<{\f{6-3p}{4p}}<\f34$ and $\f{4p}{2-p}>4$. As for the other term $r\p_r\big(\f{\psi}{r}\big)$, our strategy is to testify $r\p_r\big(\f{\psi}{r}\big)\in W_{\rm loc}^{1,p}(\R^2_{+})$, which is based on the inequality
$$\|r\p_r\big(\f{\psi}{r}\big)\|_{L_{\rm loc}^{\f{2p}{2-p}}(\R^3)}\leq C\|r\p_r\big(\f{\psi}{r}\big)\|_{L_{\rm loc}^{\f{2p}{2-p}}(\R_{+}^2)}$$
and the Sobolev embedding inequality $W_{\rm loc}^{1,p}(\R^2_{+})\hookrightarrow L_{\rm loc}^{\f{2p}{2-p}}(\R^2_{+})$ for any $1<p<2$. To start with, we recall
\eqref{uz} that $r\p_r\big(\f{\psi}{r}\big)=u_z-\f{2\psi}r$. Effectively, in Lemma \ref{estuu3}, we have proved $u_z\in L_{\rm loc}^{p}(\R^2_{+})$. Besides, the $L_{\rm loc}^{p}(\R^2_{+})$ estimate of $\f{\psi}r$ was also done in \eqref{estu35}, that can be summarized in the following estimate
\ben\label{estu46}
\|r\p_r\big(\f{\psi}{r}\big)\|_{L^p([0,R]\times[-R,R])}\leq C(R+1)\|\f{w^\theta}{r}\|_{L^1\cap L^p(\R^3)}.
\een
In next stage, to prove $\tilde{\nabla}\big[r\p_r\big(\f{\psi}{r}\big)\big]\in L_{\rm loc}^{p}(\R^2_{+})$, we will do some decompositions, which thereby make Lemma \ref{axi1} effective. More precisely, we will prove $\p_r\big(r\p_r\big(\f{\psi}{r}\big)\big)=r\p_r^2\big(\f{\psi}{r}\big)+\p_r\big(\f{\psi}{r}\big),\,\p_z\big(r\p_r\big(\f{\psi}{r}\big)\big)=r\p_{rz}^2\big(\f{\psi}{r}\big)\in L_{\rm loc}^{p}(\R^2_{+}).$ In this end, we first list the inequality
$$\|f\|_{L_{\rm loc}^p(\R_{+}^2)}\leq C\|\f{f}r\|_{L_{\rm loc}^p(\R^3)}$$
that holds for any function $f=f(r,z,t)$. This means that it suffices to verify $\f1r\p_r\big(r\p_r\big(\f{\psi}{r}\big)\big)=\p_r^2\big(\f{\psi}{r}\big)+\f1r\p_r\big(\f{\psi}{r}\big),\,\f1r\p_z\big(r\p_r\big(\f{\psi}{r}\big)\big)=\p_{rz}^2\big(\f{\psi}{r}\big)\in L_{\rm loc}^{p}(\R^3)$, which certainly holds according to Lemma \ref{axi1}. Thus, we finish all the proof.
\endproof

\vskip .1in
Thus, for $1<p<2$, we have established the ${L^{\f{2p}{2-p}}_{\rm loc}(\R^3)}$ estimates of velocity fields. When $p\geq 2$, it is well known that the Sobolev embedding $W_{\rm loc}^{1,p}(\R^3)\hookrightarrow L_{\rm loc}^6(\R^3)$ holds, which also helps us deriving the following conclusion.

\begin{lemma}\label{estuu5}Let $u=u_r(r,z,t)e_r+u_z(r,z,t)e_z$ be a smooth axisymmetric velocity fields without swirl, $\f{w_0^\theta}r\in L^1\cap L^p(\R^3)$
with some $p > 1$, then there exists an $\alpha>0$ depending only on $p$ such that $u\in L_{\rm loc}^{2+\alpha}(\R^3)$.
\end{lemma}

\vskip .3in
\section{Global existence of weak solutions}\hspace*{\parindent}

This section is devoted to the global existence of weak solutions. The first step is to construct a family of approximation solutions. To begin with, we would like to introduce the standard mollifier $\rho_{\epsilon}$, which can be described by
$$\rho_{\epsilon}(x)=\f{1}{{\epsilon}^3} \rho(\f{|x|}{\epsilon}),$$
where $\rho\in C_0^{\infty}(\R^3),\,\,\, \rho\geq0,\,\,\,{\rm supp}\,\rho\in\{|x|\leq 1\}$ and $\int_{\R^3}\rho\,dx=1.$ Then, we define a cut-off function
$\chi_{\epsilon}$ by
$$\chi_{\epsilon}(x)=\chi(\f{|x|}{\epsilon}),$$
where $\chi\in C_0^{\infty}(\R^3),\,\,\, 0\leq\chi\leq1$, and $\chi(x)=1\,\,\,{\rm on}\,\,\,\{|x|\leq 1\}$, $\chi(x)=0\,\,\,{\rm on}\,\,\,\{|x|\geq 2\}$. Through borrowing these definitions, we then drive the following theorem.

\vskip .1in
\begin{theorem}\label{exist} Given an initial data $w_0=w_0^{\theta}e_{\theta}$ such that $\f{w_0^{\theta}}r\in{L^1\cap L^p(\R^3)}$ for some $p>1,$ then there exists a family of smooth axisymmetric solutions $u^{\epsilon}$ with zero swirl component and initial data $u_0^{\epsilon}$ for any $T>0.$ Here, $w_0^{\epsilon}(x)=\rho_{\epsilon}\ast w_0(x)$ and $u_0^{\epsilon}=\nabla\times{(-\Delta)}^{-1}w_0^{\epsilon}$. In addition, it holds that
\ben\label{exist1}
\|u^{\epsilon}\|_{W^{1,p}(B_{R}\times[-R,R])}\leq C_{R}
\een
and
\ben\label{exist2}
\|{u}^{\epsilon}\|_{L^{2+\alpha}(B_R\times[-R,R])}\leq C_{R},
\een
where $\alpha$ be as in Lemma \ref{estuu5} and $C_R$ be the constants depending only on $R$.
\end{theorem}
\no{\bf Proof.}\quad Initially, we construct
$$w_0^{\epsilon}=\chi_{\epsilon}(x)({\rho_{\epsilon}}\ast w_0)(x).$$
According to our construction for initial data, it is clear that $w_0^{\epsilon}$ is axisymmetric. Then we denote by $u_0^{\epsilon}$ the corresponding velocity determined by the {\it Biot-Savart law}, namely $u_0^{\epsilon}=\nabla\times{(-\Delta)}^{-1}w_0^{\epsilon}$. Again by our assumptions on initial data,  $\nabla\times u_0^{\epsilon}=w_0^{\epsilon}$ has only swirl component $w_{\theta}^{\epsilon}(0,x)$ such that $w_0^{\epsilon}=w_{\theta}^{\epsilon}(0,x)e_\theta$. Therfore, it is clear to conclude that $u_0^{\epsilon}$ has zero swirl component, i.e., $u_{\theta}^{\epsilon}(0,x)=0$. Moreover, $u_0^{\epsilon}\in C^{\infty}(\R^3)$ and belongs to the space $V=\{u\in H^3(\R^3)|\,\nabla\cdot u=0\}.$

Subsequently, by Theorem 2.4 of \cite{DiPerna1}, there exists a unigue global axisymmetric smooth solution $u^{\epsilon}$.
What's more, because Euler equations keep invariant under the rotation and translation transformations, it is obvious
that the vector fields $u^{\epsilon}$ is still axisymmetric. Besides, the swirl component $u_\theta^{\epsilon}$ is also vanishing due to its initial data $u_{0,\theta}^{\epsilon}$ given zero.

Finally, we recall a well known conclusion that
\ben\label{winitial}
\|\f{w_0^{\epsilon}}r\|_{L^p(R^3)}\leq\|\f{\rho_{\epsilon}\ast w_0{^\theta}}r\|_{L^p(\R^3)}\leq C\|\f{w_0^{\theta}}r\|_{L^p(\R^3)},\quad \forall p\in[1,\infty],
\een
whose proof can be referred to Lemma A.1 in \cite{BenDan}. Thus, through evoking the transport equation (\ref{equ4}) satisfied by $\f{w_{\theta}^{\epsilon}}r$, applying (\ref{est1}) and \eqref{winitial}, we can conclude that $\|\f{w_{\theta}^{\epsilon}}r\|_{L^1\cap L^p(R^3)}\leq C.$  This together with Proposition \ref{estuu1}-\ref{estuu4} leads to \eqref{exist1} and \eqref{exist2}.
\endproof

\vskip .1in
As discussed in the introduction, in order to prove the main theorem, it suffices to build up the strong convergence of approximating solutions in the space $L^2(0, T; L_{\rm loc}^2(\R^3))$. Based on it, for the approximating solutions we constructed, one can then take the limit in the sense of Definition \ref{weakE}, which is essential in establishing the global existence of weak solutions. In the end, with the help of {\it a priori estimates} in Proposition \ref{estuu1}-\ref{estuu4}, we get to prove our main theorem as follow.\\

\no{\bf Proof of Theorem \ref{thm41}.}\quad As stated in the introduction, for any $p>1$, the $W_{\rm loc}^{1,p}(\R^3)$ estimates of velocity fields can not guarantee the strong convergence of approximating solutions in $L^2(0, T; L_{\rm loc}^2(\R^3))$, but in $L^{2}(0,T; Q)$ for any $Q\subset\subset {\R^3}\backslash\{x\in R^3|r = 0\}.$ Hence, we will verify the strong convergence by dividing any local domain of $\R^3$ into two parts: the region near the axis of symmetry, and the region away from it. On one hand, thanks to Lemma \ref{estuu5}, for the approximating solutions constructed in Theorem \ref{exist}, there exists $u$ such that
that
\ben\label{weakcon}
{u^{\epsilon}}\rightharpoonup u\quad {\rm in}\quad L^{\infty}(0,T;L^{2+\alpha}(B_R\times[-R,R])).
\een

On the other hand, for the region $C_R\times[-R,R]=\{(x,y)\in\R^2|\f1R\leq\sqrt{x^2+y^2}\leq R\}\times[-R,R]$, it clearly holds $\|u^{\epsilon}\|_{L^{\infty}(0, T; W^{1,p}(C_R\times[-R,R]))}\leq C_R$ by Theorem \ref{exist}. Then by using equation $(\ref{equ1})^1$, it further holds $\|\p_t u^{\epsilon}\|_{L^{\infty}(0, T; W^{-1,p^{*}}(C_R\times[-R,R]))}\leq C_R,$ where $p^{*}=\f{p}{p-1}.$ Then by noticing that $|u|$ is a function of variables $r$, $z$ and $t$, one can conlude that
$$\|u^{\epsilon}\|_{L^{\infty}(0, T; W^{1,p}([\f1R,R]\times[-R,R];drdz))}+\|\p_t u^{\epsilon}\|_{L^{\infty}(0, T; W^{-1,p^{*}}([\f1R,R]\times[-R,R];drdz))}\leq C(R).$$
Next, by applying the Aubin-Lions lemma and compact embeddings $W^{1,p}([\f1R,R]\times[-R,R])$
$\hookrightarrow L^{2}([\f1R,R]\times[-R,R])$ for any $p>1$, we can then find a subsequence $u^{\epsilon_j}$ (depending on $R$) such that
$${u^{\epsilon_j}}\rightarrow \bar{u} \quad {\rm in} \quad L^{2}(0,T; ([\f1R,R]\times[-R,R];drdz)).$$
Then, by the diagonal selection process, one can then extract a subsequence of $u^{\epsilon_j}$ independent of $R$ (still denoted by $u^{\epsilon_j}$) such that
$$\|u^{\epsilon_j}-\bar{u}\|_{L^{2}(0,T; ([\f1R,R]\times[-R,R];drdz))}\rightarrow 0\quad{\rm as}\quad {\epsilon_j\rightarrow0},$$
which also implies that
$$\|u^{\epsilon_j}-\bar{u}\|_{L^{2}(0,T; C_R\times[-R,R])}\rightarrow 0\quad{\rm as}\quad {\epsilon_j\rightarrow0}.$$
This means ${u^{\epsilon_j}}\rightarrow \bar{u}$ in $L^{2}(0,T;Q)$.
for any $Q\subset\subset {B_R\times[-R,R]}\backslash\{x\in R^3|r = 0\}.$ Then by considering the uniqueness of limits and \eqref{weakcon}, we actually have derived
\ben\label{localcon}
{u^{\epsilon_j}}\rightarrow {u}\quad{\rm in}\quad L^{2}(0,T;Q).
\een

Now, it suffices to verify the strong convergence of velocity fields in $L^{2}(0,T;B_R\times[-R,R])$.
For any $\epsilon>0$, we firstly take $Q\subset\subset {B_R\times[-R,R]}\backslash\{x\in R^3|r = 0\}$ such that the measure
$\mu(B_R\times[-R,R]\backslash Q)<\big(\f{\epsilon}{2\sqrt{2T}C_R}\big)^{\f{4+2\alpha}{\alpha}}$ for $\alpha>0$. Then according to \eqref{localcon}, there exists
$M$ such that when $j>M$, $\|u^{\epsilon_j}-{u}\|_{L^{2}(0,T; Q)}<\f{\epsilon}2$. Thus, by employing H\"{o}lder inequality and \eqref{weakcon}, for $j>M$, one further has
\beno
&&\big[\int_{0}^{T}\int_{B_R\times[-R,R]}|u^{\epsilon_j}-{u}|^2dxdt\big]^{\f1{2}}\\
&\leq&\big[\int_{0}^{T}\int_{B_R\times[-R,R]\backslash Q}|u^{\epsilon_j}-{u}|^2dxdt\big]^{\f1{2}}+\big[\int_{0}^{T}\int_{ Q}|u^{\epsilon_j}-{u}|^2dxdt\big]^{\f1{2}}\\
&\leq&\sqrt{2T}\big[\int_{B_R\times[-R,R]\backslash Q}|u^{\epsilon_j}|^2dx+\int_{B_R\times[-R,R]\backslash Q}|u|^2dx\big]^{\f1{2}}+\f{\epsilon}2\\
&\leq&\sqrt{2T}\big[\|u^{\epsilon_j}\|_{L^{2+\alpha}(B_{R}\times[-R,R])}+\|u\|_{L^{2+\alpha}(B_{R}\times[-R,R])}\big]\big[\mu(B_R\times[-R,R]\backslash Q)\big]^{\f{\alpha}{4+2\alpha}}+\f{\epsilon}2\\
&<&\epsilon.
\eeno

Until now, we actually have proved that there exists an axisymmetric vector field $u$ without swirl, such that
$$u^{\epsilon_i}\rightarrow u\quad {\rm strongly}\quad {\rm in}\quad L^2(0, T; L_{\rm loc}^2(\R^3)).$$

The last step is to pass limit in the equations \eqref{equ1} satisfied by $u^\epsilon$. As a matter of fact, it suffices to show the convergence of nonlinear term. Considering that ${u^{\epsilon_j}}\rightarrow u$ strongly in $L^2(0,T;L^2_{\rm loc}(\R^3))$, it is not hard to infer that
$$\int_0^T\int_{\R^3}u^{\epsilon_j} \cdot\nabla {\varphi}\cdot u^{\epsilon_j}\,dxdt\rightarrow \int_0^T\int_{\R^3}u \cdot\nabla {\varphi}\cdot u\,dxdt$$
for any $\varphi\in C_0^\infty((0,T];\R^3)$. This shows that $u$ is a weak solution of 3D incompressible axisymmetric Euler equations without swirl in the sense of Definition \ref{weakE}.

\endproof

\vskip .4in
\section*{Acknowledgments}
This work was started when the second author was doing his postdoctoral research supported by the CNPq grant $\sharp$ 501376/2013-1 at Federal University of Rio de Janeiro (UFRJ) of Brazil, and he would like to thank Prof. Milton C. Lopes Filho and Prof. Helena J. Nussenzveig Lopes for their hosting and hospitality. The second author also would like to thank Prof. Edriss S. Titi for his valuable suggestions about this problem when he was visiting UFRJ. Jiu is supported by National Natural Sciences Foundation of China (No. 11171229, No. 11231006). Liu is supported by the Connotation Development Funds of Beijing University of Technology. Niu is supported by National Natural Sciences Foundation of China (No. 11471220) and Beijing Municipal Commission of Education grants (No. KM201610028001).
\par

\end{document}